\documentclass[oneside,11pt]{amsart}
\usepackage{amsmath, amsfonts,amsthm,times,graphics}
\usepackage{mathrsfs}
\usepackage[usenames]{color}
\usepackage[active]{srcltx}
 \makeatletter
\renewcommand*\subjclass[2][2000]{%
  \def\@subjclass{#2}%
  \@ifundefined{subjclassname@#1}{%
    \ClassWarning{\@classname}{Unknown edition (#1) of Mathematics
      Subject Classification; using '1991'.}%
  }{%
    \@xp\let\@xp\subjclassname\csname subjclassname@#1\endcsname
  }%
}
 \makeatother

\newtheorem*{ThmA}{Theorem A}
\newtheorem*{ThmB}{Theorem B}
\newtheorem*{ThmC}{Theorem C}

\newtheorem*{ThmE}{Theorem E}

\newtheorem*{LemD}{Lemma D}

\newtheorem{Thm}{Theorem}[section]
\newtheorem{Cor}[Thm]{Corollary}
\newtheorem{Lem}[Thm]{Lemma}
\newtheorem{Pro}[Thm]{Proposition}
\theoremstyle{definition}

\theoremstyle{remark}

\numberwithin{equation}{section}

\newcommand{\ee}{\mathrm{e}}

\theoremstyle{definition}

\def\be{\begin{equation}}
\def\ee{\end{equation}}

\newcommand{\ben}{\begin{enumerate}}
\newcommand{\een}{\end{enumerate}}

\newcommand{\br}{\begin{rem}}
\newcommand{\er}{\end{rem}}
\newcommand{\brs}{\begin{rems}}
\newcommand{\ers}{\end{rems}}
\newcommand{\bo}{\begin{obser}}
\newcommand{\eo}{\end{obser}}
\newcommand{\bos}{\begin{obsers}}
\newcommand{\eos}{\end{obsers}}
\newcommand{\bpf}{\begin{pf}}
\newcommand{\epf}{\end{pf}}
\newcommand{\ba}{\begin{array}}
\newcommand{\ea}{\end{array}}
\newcommand{\beq}{\begin{eqnarray}}
\newcommand{\beqq}{\begin{eqnarray*}}
\newcommand{\eeq}{\end{eqnarray}}
\newcommand{\eeqq}{\end{eqnarray*}}

\numberwithin{equation}{section}

\newcounter{minutes}
\divide\time by 60
\newcounter{hours}
\multiply\time by 60 \addtocounter{minutes}{-\time}
\begin{document}
\title{Characterizations of composition operators on Bloch and Hardy type spaces}

\author[Shaolin Chen and  Hidetaka Hamada]{Shaolin Chen and Hidetaka Hamada}

\address{S. L.  Chen, College of Mathematics and
Statistics, Hengyang Normal University, Hengyang, Hunan 421002,
People's Republic of China; Hunan Provincial Key Laboratory of
Intelligent Information Processing and Application,  421002,
People's Republic of China.} \email{mathechen@126.com}

\address{H. Hamada, Faculty of Science and Engineering, Kyushu Sangyo University,
3-1 Matsukadai 2-Chome, Higashi-ku, Fukuoka 813-8503, Japan.}
\email{ h.hamada@ip.kyusan-u.ac.jp}


\maketitle

\def\thefootnote{}
\makeatletter\def\thefootnote{\@arabic\c@footnote}\makeatother


\begin{abstract}
The main purpose of this paper is to investigate
characterizations of composition operators on Bloch and Hardy type spaces. Initially,
we  use general doubling weights to study
the composition operators from
 harmonic Bloch type spaces on the unit disc $\mathbb{D}$
 to pluriharmonic Hardy spaces on the Euclidean unit ball $\mathbb{B}^n$.  Furthermore,
we develop some new methods to study the composition operators from harmonic Bloch type spaces on $\mathbb{D}$
to pluriharmonic Bloch type spaces on $\mathbb{D}$.
 Additionally, some application to new characterizations of the composition operators
between pluriharmonic Lipschitz type spaces to be bounded or compact will be presented. The obtained results of this paper provide the improvements and
extensions of the corresponding known results.
\end{abstract}

\maketitle \pagestyle{myheadings} \markboth{ S. L. Chen and H.
Hamada}{Characterizations of composition operators on Bloch and Hardy type spaces}

\maketitle

\section{ Introduction}\label{sec1}

The study of composition operators on various Banach spaces of holomorphic
functions and planar harmonic functions is currently a very active field of complex and functional analysis (see \cite{AD,CHZ2022MZ,CPR,GZ,HO,Kw,MM,PX,Sha,Z2}).
 This paper continues the study of previous work of authors \cite{CHZ2022MZ,CPR} and
 is mainly motivated by the articles of Kwon \cite{Kw},  Pavlovi\'c \cite{Pav-2008},  Wulan et al. \cite{WZZ} and Zhao \cite{Zhao}.
First,
  we  use general doubling weights to study
the composition operators from
 harmonic Bloch type spaces on the unit disc $\mathbb{D}$
 to pluriharmonic Hardy spaces on the Euclidean unit ball $\mathbb{B}^n$.  In addition,
we develop some new methods to study the composition operators from harmonic Bloch type spaces on $\mathbb{D}$
to pluriharmonic Bloch type spaces on $\mathbb{D}$.
 At last, some application to new characterizations of the composition operators
between pluriharmonic Lipschitz type spaces to be bounded or compact will be given. The obtained results of this paper provide the improvements and
extensions of the corresponding known results.
 In particular, we improve and extend the main results of Chen et al. \cite{CHZ2022MZ} and Kwon \cite{Kw}, and we also establish a completely different characterization from  Pavlovi\'c \cite{Pav-2008}.
  In order to state our main results, we need to recall some basic definitions and introduce some necessary terminologies.

Let $\mathbb{C}^{n}$  be the complex space of dimension $n$,  and let $\mathbb{C}:=\mathbb{C}^{1}$ be the complex plane,
 where $n$ is a positive integer. For $z=(z_{1},\ldots,z_{n})\in\mathbb{C}^{n}$ and $w=(w_{1},\ldots,w_{n})\in\mathbb{C}^{n} $, we
write $\langle z,w\rangle := \sum_{k=1}^nz_k\overline{w}_k$ and $
|z|:={\langle z,z\rangle}^{1/2}.$ For $a\in
\mathbb{C}^n$, we set $\mathbb{B}^n(a, r)=\{z\in \mathbb{C}^{n}:\,
|z-a|<r\}. $ In particular, let $\mathbb{B}^n:=\mathbb{B}^n(0, 1)$
and $\mathbb{D}:=\mathbb{B}^1$.

A twice continuously differentiable complex-valued function $f$
defined on a domain $\Omega\subset\mathbb{C}^{n}$ is called a
pluriharmonic function if for each fixed $z\in \Omega$ and
$\theta\in\partial\mathbb{B}^{n}$, the function
$f(z+\zeta\theta)$ is harmonic in $\{\zeta:\; |\zeta|<
d_{\Omega}(z)\}$,
where $d_{\Omega}(z)$ is the distance from
$z$ to the boundary $\partial\Omega$ of $\Omega$.
Recently, pluriharmonic functions have been widely studied (cf.
\cite{CH-2022,CHPV,DHK-2011,Iz,Ra,RU,Ru-1,Vl}).
If
$\Omega\subset\mathbb{C}^{n}$ is a simply connected domain,
then a function $f:\,\Omega\rightarrow
\mathbb{C}$ is pluriharmonic if and only if $f$ has a decomposition
$f=h+\overline{g},$ where $h$ and $g$ are holomorphic in $\Omega$
 (see \cite{Vl}). This decomposition is unique up to an additive constant.
 From this decomposition, it is easy to know that the class of pluriharmonic
  functions is broader than that of holomorphic functions. Furthermore, a twice continuously
differentiable real-valued function in a simply connected domain
$\Omega$ is pluriharmonic if and only if it is the real part of some
holomorphic function on $\Omega$.
Obviously, all pluriharmonic functions are harmonic.
In particular, if $n=1$, then
the converse holds  (cf. \cite{Du}). Throughout of this paper, we use $\mathscr{H}(\Omega)$ and
$\mathscr{PH}(\Omega)$ to denote the set of all holomorphic
functions of a domain $\Omega\subset\mathbb{C}^{n}$ into $\mathbb{C}$ and that of all pluriharmonic
functions of  $\Omega$  into $\mathbb{C}$, respectively.

\subsection{Hardy spaces}

For $p\in(0,\infty]$, the pluriharmonic Hardy space $\mathscr{PH}^{p}(\mathbb{B}^{n})$ consists of all those
functions $f\in\mathscr{PH}(\mathbb{B}^{n})$ such that, for $p\in(0,\infty)$,
$$\|f\|_{p}:=\sup_{r\in[0,1)}M_{p}(r,f)<\infty,$$ and, for $p=\infty$, $$\|f\|_{\infty}:=\sup_{r\in[0,1)}M_{\infty}(r,f)<\infty,$$
where
$$M_{p}(r,f)=\left(\int_{\partial\mathbb{B}^{n}}|f(r\zeta)|^{p}d\sigma(\zeta)\right)^{\frac{1}{p}},
~M_{\infty}(r,f)=\sup_{\zeta\in\partial\mathbb{B}^{n}}|f(r\zeta)|$$
and $d\sigma$ denotes the normalized Lebesgue surface measure on
$\partial\mathbb{B}^{n}$. In particular, we use
$\mathscr{H}^{p}(\mathbb{B}^{n}):=\mathscr{H}(\mathbb{B}^{n})\cap\mathscr{PH}^{p}(\mathbb{B}^{n})$
to denote the holomorphic Hardy space. If
$f\in\mathscr{PH}^{p}(\mathbb{B}^{n})$ for some $p\in (1,\infty)$,
then the radial limits
$$f(\zeta)=\lim_{r\rightarrow1^{-}}f(r\zeta)$$ exist for almost every $\zeta\in\partial\mathbb{B}^{n}$ (see \cite[Theorems 6.7,  6.13 and 6.39]{ABR-2001}).
Moreover, if $p\in[1,\infty)$, then $\mathscr{PH}^{p}(\mathbb{B}^{n})$ is a normed space
with respect to the norm $\| \cdot\|_p$
and if $p\in(0,1)$,
then
$\mathscr{PH}^{p}(\mathbb{B}^{n})$ is a metric space with respect to the metric $d(f,g)=\| f-g\|_p^p$.

\subsection{Bloch type spaces}
For a
pluriharmonic function $f$ of $\mathbb{B}^{n}$ into $\mathbb{C}$, let
$$\nabla\,f:=\left(\frac{\partial\,f}{\partial\,z_{1}},\ldots,\frac{\partial\,f}{\partial\,z_{n}}\right),~ \overline{\nabla}\,
f:=\left(\frac{\partial\,f}{\partial\,\overline{z}_{1}},\ldots,\frac{\partial\,f}{\partial\,\overline{z}_{n}}\right)$$
and
$$\Lambda_{f}(z):=\max_{\theta\in\partial\mathbb{B}^{n}}|\langle\nabla\,f(z),\overline{\theta}\rangle+\langle\overline{\nabla}\,f(z),\theta\rangle|$$
 for $z=(z_{1},\ldots,z_{n})\in\mathbb{B}^{n}$.

 A continuous non-decreasing function $\omega:~[0,1)\rightarrow(0,\infty)$ is called a
 weight if $\omega$ is  unbounded (see \cite{AD}). Moreover,
 a weight  $\omega$ is called doubling if there is a constant $C>1$
such that
$$\omega(1-s/2)<C\omega(1-s)$$
for $s\in(0,1]$.

For a weight $\omega$, we use $\mathscr{B}_{\omega}(\mathbb{B}^n)$ to denote
the pluriharmonic Bloch type space consisting of all
complex-valued pluriharmonic functions defined in $\mathbb{B}^n$ with the norm
$$\|f\|_{\mathscr{B}_{\omega}(\mathbb{B}^n)}:=
|f(0)|+\sup_{z\in \mathbb{B}^n}\mathscr{B}^{f}_{\omega}(z)<\infty,$$
where  $\mathscr{B}^{f}_{\omega}(z)=\Lambda_{f}(z)/\omega(|z|).$ It is easy to know that $\mathscr{B}_{\omega}(\mathbb{B}^n)$
is a complex Banach space.   Furthermore, let
$$\|f\|_{\mathscr{B}_{\omega}(\mathbb{B}^n),s}:=
\sup_{z\in\mathbb{B}^n}\mathscr{B}^{f}_{\omega}(z)<\infty$$ be the semi-norm. If $f\in \mathscr{B}_{\omega}(\mathbb{B}^n)$, then we call
$f$ a pluriharmonic Bloch function.

\subsection{Composition operators}
Given a holomorphic  function $\phi$ of $\mathbb{B}^{n}$ into $\mathbb{D}$, the composition
operator $C_{\phi}:~\mathscr{PH}(\mathbb{D})\rightarrow\mathscr{PH}(\mathbb{B}^{n})$ is defined by
$$C_{\phi}(f)=f\circ\phi,$$ where $f\in\mathscr{PH}(\mathbb{D})$.

In 1987, Shapiro \cite{Sha} gave a complete characterization of
compact composition operators on $\mathscr{H}^{2}(\mathbb{D})$, with
a number of interesting consequences for peak sets, essential norm
of composition operators, and so on. Recently, the studies of
composition operators on holomorphic function spaces have been
attracted much attention of many mathematicians (see
\cite{AD,CHZ2022MZ,CPR,GZ,HO,Kw,MM,PX,Z2}).  In particular, Kwon \cite{Kw}
investigated some characterizations of composition operators from the holomorphic
Bloch spaces to the holomorphic Hardy spaces
to be bounded or compact.
Let us recall the main result in
\cite{Kw} as follows.

For $\zeta\in\partial\mathbb{B}^{n}$ and $\alpha\in(1,\infty)$, we use $D_{\alpha}^{n}(\zeta)$
 to denote the Koranyi approach domain defined by
$$D_{\alpha}^{n}(\zeta)=\left\{z\in\mathbb{B}^{n}:~|1-\langle\,z,\zeta\rangle|<\frac{\alpha}{2}(1-|z|^{2})\right\}.$$
For $\phi:~\mathbb{B}^{n}\rightarrow\mathbb{D}$ holomorphic,
let $$M_{\alpha}\phi(\zeta)=\sup\left\{\log\frac{1}{1-|\phi(z)|^{2}}:~z\in\,D_{\alpha}^{n}(\zeta)\right\}$$ be the maximal function (see \cite{Ko}).
Let $\mbox{Aut}(\mathbb{B}^n)$ denote the group of holomorphic automorphisms of $\mathbb{B}^n$.

\begin{ThmA}{\rm (\cite[Theorems 5.1 and 5.10]{Kw})}\label{Kwon-2002}
Let $p\in(0,\infty)$, $1<\alpha,~\beta<\infty$ and $\omega(t)=1/(1-t^{2})$ for $t\in[0,1)$. If $\phi:~\mathbb{B}^{n}\rightarrow\mathbb{D}$ is holomorphic, then
the followings are equivalent:
\begin{enumerate}
\item[{\rm (1)}] $\sup_{r\in(0,1)}\int_{\partial\mathbb{B}^{n}}\left(\frac{1}{2}\log\frac{1+|\phi(r\zeta)|}{1-|\phi(r\zeta)|}\right)^{\frac{p}{2}}d\sigma(\zeta)<\infty$;

\item[{\rm (2)}] $\int_{\partial\mathbb{B}^{n}}\left(M_{\beta}\phi(\zeta)\right)^{\frac{p}{2}}d\sigma(\zeta)<\infty$;

\item[{\rm (3)}] $\int_{\partial\mathbb{B}^{n}}\left(\int_{0}^{1}\frac{|\nabla(\phi\circ\varphi_{r\zeta})(0)|^{2}}{(1-|\phi(r\zeta)|^{2})^{2}}
    \frac{dr}{1-r}\right)^{\frac{p}{2}}d\sigma(\zeta)<\infty,$ where $\varphi_{r\zeta}\in\mbox{Aut}(\mathbb{B}^{n})$ with $\varphi_{r\zeta}(0)=r\zeta$;

    \item[{\rm (4)}] $\int_{\partial\mathbb{B}^{n}}\left(\int_{D_{\alpha}^{n}(\zeta)}
    \frac{|\nabla(\phi\circ\varphi_{z})(0)|^{2}}{(1-|\phi(z)|^{2})^{2}} \frac{dV(z)}{(1-|z|^{2})^{n+1}}\right)^{\frac{p}{2}}d\sigma(\zeta)<\infty,$
    where $dV$ denotes the  Lebesgue volume measure of $\mathbb{C}^{n}$, and $\varphi_{z}\in\mbox{Aut}(\mathbb{B}^{n})$ with $\varphi_{z}(0)=z$;

\item[{\rm (5)}] $C_{\phi}:~\mathscr{B}_{\omega}(\mathbb{D})\cap \mathscr{H}(\mathbb{D})\rightarrow\mathscr{H}^{p}(\mathbb{B}^{n})$ is a bounded operator;

\item[{\rm (6)}] $C_{\phi}:~\mathscr{B}_{\omega}(\mathbb{D})\cap \mathscr{H}(\mathbb{D})\rightarrow\mathscr{H}^{p}(\mathbb{B}^{n})$ is compact.

\end{enumerate}
\end{ThmA}

Also,
the study of composition operators and linear operators on holomorphic function spaces by using weight has
aroused great interest of many mathematicians (cf. \cite{CHZ2022MZ,GZ,Mo-20,PR-2013}).
However, there are few literatures on the theory of composition
operators of  harmonic functions. In the following, by using general doubling
 weights, we will establish the characterizations of
composition operators from harmonic Bloch type spaces
on $\mathbb{D}$ to pluriharmonic Hardy spaces on $\mathbb{B}^n$
to be bounded or compact.

\begin{Thm}\label{thm-2}
Let $p\in(0,\infty)$, $\alpha\in(1,\infty)$ and $\omega$ be a doubling function. If $\phi:~\mathbb{B}^{n}\rightarrow\mathbb{D}$ is holomorphic,
 then the followings are equivalent:
\begin{enumerate}
\item[{\rm (1)}] $C_{\phi}:~\mathscr{B}_{\omega}(\mathbb{D})\rightarrow\mathscr{PH}^{p}(\mathbb{B}^{n})$ is a bounded operator;

\item[{\rm (2)}] $ \int_{\partial\mathbb{B}^{n}}\left(\int_{0}^{1}|\nabla\phi(r\zeta)|^{2}\omega^{2}(|\phi(r\zeta)|)(1-r)
dr\right)^{\frac{p}{2}}\,d\sigma(\zeta)<\infty;
$

\item[{\rm (3)}] $ \int_{\partial\mathbb{B}^{n}}\left(\int_{D_{\alpha}^{n}(\zeta)}|\nabla\phi(z)|^{2}\omega^{2}(|\phi(z)|)(1-|z|)^{1-n}
dV(z)\right)^{\frac{p}{2}}\,d\sigma(\zeta)<\infty;$

\item[{\rm (4)}] $C_{\phi}:~\mathscr{B}_{\omega}(\mathbb{D})\rightarrow\mathscr{PH}^{p}(\mathbb{B}^{n})$ is compact.
\end{enumerate}
\end{Thm}

If we take $\omega(t)=1/(1-t^{2})$ for $t\in[0,1)$ in Theorem
\ref{thm-2}, then we extend Theorem A 
into the
following form.
\begin{Cor}
Let $p\in(0,\infty)$ and $1<\alpha,~\beta<\infty$. If
$\phi:~\mathbb{B}^{n}\rightarrow\mathbb{D}$ is holomorphic, then the
followings are equivalent:
\begin{enumerate}
\item[{\rm (1)}] $\sup_{r\in(0,1)}\int_{\partial\mathbb{B}^{n}}\left(\frac{1}{2}\log\frac{1+|\phi(r\zeta)|}{1-|\phi(r\zeta)|}\right)^{\frac{p}{2}}d\sigma(\zeta)<\infty$;

\item[{\rm (2)}] $\int_{\partial\mathbb{B}^{n}}\left(M_{\beta}\phi(\zeta)\right)^{\frac{p}{2}}d\sigma(\zeta)<\infty$;

\item[{\rm (3)}] $\int_{\partial\mathbb{B}^{n}}\left(\int_{0}^{1}\frac{|\nabla(\phi\circ\varphi_{r\zeta})(0)|^{2}}{(1-|\phi(r\zeta)|^{2})^{2}}
    \frac{dr}{1-r}\right)^{\frac{p}{2}}d\sigma(\zeta)<\infty,$ where $\varphi_{r\zeta}\in\mbox{Aut}(\mathbb{B}^{n})$ with $\varphi_{r\zeta}(0)=r\zeta$;

    \item[{\rm (4)}] $\int_{\partial\mathbb{B}^{n}}\left(\int_{D_{\alpha}^{n}(\zeta)}
    \frac{|\nabla(\phi\circ\varphi_{z})(0)|^{2}}{(1-|\phi(z)|^{2})^{2}} \frac{dV(z)}{(1-|z|^{2})^{n+1}}\right)^{\frac{p}{2}}d\sigma(\zeta)<\infty,$
    where  $\varphi_{z}\in\mbox{Aut}(\mathbb{B}^{n})$ with $\varphi_{z}(0)=z$;

\item[{\rm (5)}] $ \int_{\partial\mathbb{B}^{n}}\left(\int_{0}^{1}|\nabla\phi(r\zeta)|^{2}\frac{(1-r)}{(1-|\phi(r\zeta)|^{2})^{2}}
dr\right)^{\frac{p}{2}}\,d\sigma(\zeta)<\infty; $

\item[{\rm (6)}] $ \int_{\partial\mathbb{B}^{n}}\left(\int_{D_{\alpha}^{n}(\zeta)}\frac{|\nabla\phi(z)|^{2}}{(1-|\phi(z)|^{2})^{2}}(1-|z|)^{1-n}
dV(z)\right)^{\frac{p}{2}}\,d\sigma(\zeta)<\infty;
$

\item[{\rm (7)}] $C_{\phi}:~\mathscr{B}_{\omega}(\mathbb{D})\rightarrow\mathscr{PH}^{p}(\mathbb{B}^{n})$ is a bounded operator;

\item[{\rm (8)}] $C_{\phi}:~\mathscr{B}_{\omega}(\mathbb{D})\rightarrow\mathscr{PH}^{p}(\mathbb{B}^{n})$ is compact.

\end{enumerate}
\end{Cor}

In the case $\phi$ maps $\mathbb{D}$ into $\mathbb{D}$,
Chen et al. \cite{CPR} proved the following result.

\begin{ThmB}{\rm (\cite[Theorem 6]{CPR})}\label{Thm-CPR}
Let  $\omega(t)=1/\left((1-t^{2})^{\alpha}\left(\log\frac{e}{1-t^{2}}\right)^{\beta}\right)$ and
$\phi:\,\mathbb{D}\rightarrow\mathbb{D}$  be an analytic function, where $\alpha\in(0,\infty)$ and $\beta\leq\alpha$.
Then the followings are equivalent:
\begin{enumerate}
\item[{\rm (1)}] $C_{\phi}:~\mathscr{B}_{\omega}(\mathbb{D})\cap \mathscr{H}(\mathbb{D})\rightarrow\mathscr{PH}^{2}(\mathbb{D})$ is a bounded operator;

\item[{\rm (2)}] $\displaystyle \frac{1}{2\pi}\int_{0}^{2\pi}\int_{0}^{1}
\frac{|\phi'(re^{i\theta})|^{2}}{(1-|\phi(re^{i\theta})|)^{2\alpha}\left(\log\frac{e}{1-|\phi(re^{i\theta})|}\right)^{2\beta}}
(1-r)\,dr\,d\theta<\infty. $
\end{enumerate}
\end{ThmB}

A continuous increasing function
$\psi~:[0,\infty)\rightarrow[0,\infty)$ with
$\psi(0)=0$ is called a majorant if $\psi(t)/t$ is
non-increasing for $t>0$ (see \cite{Dy1,P}). By
using a special doubling function
$\omega(t)=1/\psi\left((1-t^{2})^{\alpha}\left(\log\frac{e}{1-t^{2}}\right)^{\beta}\right)$
for $t\in[0,1)$, the characterization of composition
operators from $\mathscr{B}_{\omega}(\mathbb{D})$ to $\mathscr{PH}^{p}(\mathbb{D})$ to be bounded or compact
was established in \cite{CHZ2022MZ} as
follows, which is the improvement of Theorem B,
where $p\in(0,\infty)$, $\alpha\in(0,\infty)$ and $\beta\in(-\infty,\alpha]$
are constants.

\begin{ThmC}{\rm (\cite[Theorem 2.4]{CHZ2022MZ})}\label{2022-C-H}
Let $p\in(0,\infty)$, $\alpha\in(0,\infty)$,
$\beta\in(-\infty,\alpha]$ and
$\omega(t)=1/\psi\left((1-t^{2})^{\alpha}\left(\log\frac{e}{1-t^{2}}\right)^{\beta}\right)$
for $t\in[0,1)$, where $\psi$ is a majorant. If
$\phi:~\mathbb{D}\rightarrow\mathbb{D}$ is a holomorphic function,
then the followings are equivalent:

\begin{enumerate}
\item[{\rm (1)}] $C_{\phi}:~\mathscr{B}_{\omega}(\mathbb{D})\rightarrow\mathscr{PH}^{p}(\mathbb{D})$ is a bounded operator;

\item[{\rm (2)}] $\int_{0}^{2\pi}\left(\int_{0}^{1}|\phi'(re^{i\theta})|^{2}\omega^{2}(|\phi(re^{i\theta})|)(1-r)
dr\right)^{\frac{p}{2}}\frac{d\theta}{2\pi}<\infty;$

\item[{\rm (3)}] $C_{\phi}:~\mathscr{B}_{\omega}(\mathbb{D})\rightarrow\mathscr{PH}^{p}(\mathbb{D})$ is
compact.
\end{enumerate}
\end{ThmC}

In the following, by using a general doubling weight,  we  give the
characterizations of composition operators from harmonic Bloch type
spaces on $\mathbb{D}$ to harmonic Hardy spaces on $\mathbb{D}$
to be bounded or compact, which is an improvement
of Theorem C
Also, the characterizations (3), (4) and (5) are new.


\begin{Thm}\label{thm-3}
Let $p\in(0,\infty)$,
$\omega$ be a doubling function and $\phi:~\mathbb{D}\rightarrow\mathbb{D}$ be a holomorphic function. Then the followings are equivalent:
\begin{enumerate}
\item[{\rm (1)}] $C_{\phi}:~\mathscr{B}_{\omega}(\mathbb{D})\rightarrow\mathscr{PH}^{p}(\mathbb{D})$ is a bounded operator;

\item[{\rm (2)}] $\int_{0}^{2\pi}\left(\int_{0}^{1}|\phi'(re^{i\theta})|^{2}\omega^{2}(|\phi(re^{i\theta})|)(1-r)
dr\right)^{\frac{p}{2}}\frac{d\theta}{2\pi}<\infty;$

\item[{\rm (3)}] $\int_{0}^{2\pi}\left(\sum_{k=0}^{\infty}2^{-2k}|\phi'(r_{k}e^{i\theta})|^{2}
    \omega^{2}(|\phi(r_{k}e^{i\theta})|)\right)^{\frac{p}{2}}\frac{d\theta}{2\pi}<\infty,$ where $r_{k}=1-2^{-k}$;

    \item[{\rm (4)}] $\int_{0}^{2\pi}\left(\int_{0}^{1}(1-r)\sup_{0<\rho<r}
    \left(|\phi'(\rho\,e^{i\theta})|^{2}\omega^{2}(|\phi(\rho\,e^{i\theta})|)\right)
dr\right)^{\frac{p}{2}}\frac{d\theta}{2\pi}<\infty;$

\item[{\rm (5)}] $ \int_{0}^{2\pi}\left(\int_{D_{\alpha}^{1}(\zeta)}|\nabla\phi(z)|^{2}\omega^{2}(|\phi(z)|)
dA(z)\right)^{\frac{p}{2}}\,\frac{d\theta}{2\pi}<\infty,$ where $dA$ denotes the  Lebesgue area measure of $\mathbb{C}$;

\item[{\rm (6)}] $C_{\phi}:~\mathscr{B}_{\omega}(\mathbb{D})\rightarrow\mathscr{PH}^{p}(\mathbb{D})$ is
compact.
\end{enumerate}
\end{Thm}

Recently, the studies of composition operators between the classical analytic Bloch spaces have attracted much attention of many mathematicians
 (cf. \cite{MM,Mo-99,Mo-20,WZZ,Zhao}). In particular, Zhao \cite{Zhao} gave  characterizations of composition operators from the analytic $\alpha$-Bloch space
 $\mathscr{B}_{\omega_{1}}(\mathbb{D})\cap \mathscr{H}(\mathbb{D})$ to the analytic $\beta$-Bloch space $\mathscr{B}_{\omega_{2}}(\mathbb{D})\cap \mathscr{H}(\mathbb{D})$ to be bounded or compact,
 where $\alpha,~\beta\in(0,\infty)$, $\omega_{1}(t)=1/(1-t^2)^{\alpha}$ and $\omega_{2}(t)=1/(1-t^2)^{\beta}$ for $t\in[0,1)$. Pavlovi\'c \cite{Pav-2008} gave some
 derivative-free characterizations of bounded composition operators between analytic Lipschitz spaces.
 In the following, we will develop some  new methods to give new
characterizations
of the composition operators between the Bloch type spaces with weights
to be bounded or compact,
  and give some application to new characterizations of the composition operators between the Lipschitz spaces  to be bounded or compact.

For $k\in\{1,2,\ldots\}$ and a weight $\omega$, let $$\mathscr{E}_{\omega}(k)=\left\{\varrho_{k}:~\mu_{\omega,k}(\varrho_{k})=\max_{x\in[0,1)}\mu_{\omega,k}(x)\right\},$$
where  $\mu_{\omega,k}(x)=x^{k-1}/\omega(x),~x\in[0,1)$.

\begin{Pro}\label{p-1}
Let $r_1=0\in \mathscr{E}_{\omega}(1)$
and let $r_k \in \mathscr{E}_{\omega}(k)$, $k\in\{2,3,\ldots\}$, be arbitrarily chosen.
Then $\{ r_{k}\}$ is  a non-decreasing sequence.
\end{Pro}

\begin{Thm}\label{Thm-Bloch}
For $k\in\{1,2,\ldots\}$,
suppose that $\omega_{1}$ is a weight so that a point $r_{k}$ can be selected from each set $\mathscr{E}_{\omega_{1}}(k)$
with $r_1=0$ and
$$ \lim_{k\rightarrow\infty}r_{k}^{k-1}=\gamma>0.$$
Let $\omega_{2}$ be a weight, and let $\phi$ be a holomorphic function of $\mathbb{B}^n$ into $\mathbb{D}$. Then
\begin{enumerate}
\item[{\rm (1)}]
$C_{\phi}:~\mathscr{B}_{\omega_{1}}(\mathbb{D})\rightarrow\mathscr{B}_{\omega_{2}}(\mathbb{B}^n)$ is bounded if and only if
\be\label{eqx-4}\sup_{k\geq1}\frac{\omega_{1}(r_{k})}{k}\|\phi^{k}\|_{\mathscr{B}_{\omega_{2}}(\mathbb{B}^n)}<\infty.\ee
\item[{\rm (2)}]
$C_{\phi}:~\mathscr{B}_{\omega_{1}}(\mathbb{D})\rightarrow\mathscr{B}_{\omega_{2}}(\mathbb{B}^n)$ is compact if and only if
\be\label{composition-Bloch-compact}
\lim_{k\to \infty}\frac{\omega_{1}(r_{k})}{k}\|\phi^{k}\|_{\mathscr{B}_{\omega_{2}}(\mathbb{B}^n)}=0.\ee
\end{enumerate}
\end{Thm}

For $\alpha\in(0,1]$, the Lipschitz space $\mathscr{L}_{\alpha}(\mathbb{B}^n)$ consists of all those
functions $f\in\mathcal{C}(\mathbb{B}^n)$ satisfying
$$\|f\|_{\mathscr{L}_{\alpha}(\mathbb{B}^n),s}:=\sup_{z,w\in\mathbb{B}^n,z\neq\,w}\frac{|f(z)-f(w)|}{|z-w|^{\alpha}}<\infty,$$
where $\mathcal{C}(\mathbb{B}^n)$ is a set of all continuous functions defined in $\mathbb{B}^n$.


By using Theorem \ref{Thm-Bloch}, we give a characterization of the composition operators
between pluriharmonic Lipschitz type spaces to be bounded or compact
which is completely different from  Pavlovi\'c \cite{Pav-2008}  as follows.

\begin{Thm}\label{Thm-4}Suppose that $k\in\{1,2,\ldots\}$,
$\alpha,\beta\in(0,1)$ and  $\omega(t)=1/(1-t)^{1-\beta}$ for $t\in[0,1)$.
Let $\phi$ be a holomorphic function of  $\mathbb{B}^n$ into $\mathbb{D}$.
Then,
\begin{enumerate}
\item[{\rm (1)}]
$C_{\phi}:~\mathscr{L}_{\alpha}(\mathbb{D})\cap\mathscr{PH}(\mathbb{D})\rightarrow\mathscr{L}_{\beta}(\mathbb{B}^n)\cap\mathscr{PH}(\mathbb{B}^n)$
 is bounded if and only if $$\sup_{k\geq1}\left\{k^{-\alpha}\|\phi^{k}\|_{\mathscr{B}_{\omega}(\mathbb{B}^n)}\right\}<\infty;$$
\item[{\rm (2)}]
$C_{\phi}:~\mathscr{L}_{\alpha}(\mathbb{D})\cap\mathscr{PH}(\mathbb{D})\rightarrow\mathscr{L}_{\beta}(\mathbb{B}^n)\cap\mathscr{PH}(\mathbb{B}^n)$
 is compact if and only if $$\lim_{k\to \infty}\left\{k^{-\alpha}\|\phi^{k}\|_{\mathscr{B}_{\omega}(\mathbb{B}^n)}\right\}=0.$$
\end{enumerate}
\end{Thm}

The proofs of Theorems \ref{thm-2}$\sim$\ref{Thm-4} and Proposition \ref{p-1}
will be presented in Sect. \ref{sec2}.

\section{The proofs of the main results}\label{sec2}


Denote by $L^{p}(\partial\mathbb{B}^{n})$ $(p\in(0,\infty ))$ the
set of all measurable functions $F$ of $\partial\mathbb{B}^{n}$ into
$\mathbb{C}$ with

$$\|F\|_{L^{p}}=\left(\int_{\partial\mathbb{B}^{n}}|F(\zeta)|^{p}d\sigma(\zeta)\right)^{\frac{1}{p}}<\infty.$$
Given $f\in\mathscr{H}^{p}(\mathbb{B}^{n})$, the Littlewood-Paley
type $\mathscr{G}$-function is defined as follows

$$\mathscr{G}(f)(\zeta)=\left(\int_{0}^{1}|\nabla\,f(r\zeta)|^{2}(1-r)dr\right)^{\frac{1}{2}},~\zeta\in\partial\mathbb{B}^{n}.$$
Then

\be\label{g-1} f\in\mathscr{H}^{p}(\mathbb{B}^{n})~\mbox{if and only
if}~ \mathscr{G}(f)\in\,L^{p}(\partial\mathbb{B}^{n})\ee for
$p\in(0,\infty)$ (see  \cite{AB,KL,Stei}). The conclusion of (\ref{g-1})
also can be rewritten in the following form. There exists a
positive constant $C$, depending only on $p$, such that

\be\label{g-1-1}\frac{1}{C}\|f\|_{p}^{p}\leq|f(0)|^{p}+\int_{\partial\mathbb{B}^{n}}\big(\mathscr{G}(f)(\zeta)\big)^{p}d\sigma(\zeta)\leq\,C\|f\|_{p}^{p}\ee
for $p\in(0,\infty)$ (see \cite{AB,KL,Stei}). For $\alpha\in(1,\infty)$ and $p\in(0,\infty)$,  it follows
from  \cite[Theorem 3.1]{AB} (or \cite[Theorem 1.1]{KL}) that there is a positive
constant $C$ such that
\beq\label{g-1-1.1}\frac{1}{C}\int_{\partial\mathbb{B}^{n}}\big(\mathscr{G}(f)(\zeta)\big)^{p}d\sigma(\zeta)&\leq&
\int_{\partial\mathbb{B}^{n}}\left(\mathscr{A}_{\alpha}f(\zeta)\right)^{p}d\sigma(\zeta)\\ \nonumber
&\leq&\,
C\int_{\partial\mathbb{B}^{n}}\big(\mathscr{G}(f)(\zeta)\big)^{p}d\sigma(\zeta),\eeq
 where $f\in\mathscr{H}^{p}(\mathbb{B}^{n})$ and $$\mathscr{A}_{\alpha}f(\zeta)=\left(\int_{D_{\alpha}^{n}(\zeta)}|\nabla\,f(z)|^{2}(1-|z|)^{1-n}
dV(z)\right)^{1/2}.$$
It follows from (\ref{g-1-1})  and (\ref{g-1-1.1}) that there is a positive
constant $C$ such that

\be\label{g-1-1.2}\frac{1}{C}\|f\|_{p}^{p}\leq|f(0)|^{p}+\int_{\partial\mathbb{B}^{n}}\left(\mathscr{A}_{\alpha}f(\zeta)\right)^{p}d\sigma(\zeta)\leq\,C\|f\|_{p}^{p}.\ee





The following result easily follows from \cite[Lemma 1]{AD} and
\cite[Theorem 2]{AD}.

\begin{Lem}\label{lem-1}
Let  $\omega$ be a doubling function. Then there exist  functions
$f_{j}\in\mathscr{B}_{\omega}(\mathbb{D})\cap \mathscr{H}(\mathbb{D})$
$(j\in\{1,2\})$ such that, for $z\in\mathbb{D}$,
\[
\sum_{j=1}^{2}|f_{j}'(z)|\geq\omega(|z|).
\]
\end{Lem}

The following result is well-known.

\begin{LemD}{\rm (cf. \cite[Lemma 5]{CPR})}\label{Lemx}
Suppose that $a,~b\in[0,\infty)$ and $q\in(0,\infty)$. Then
$$(a+b)^{q}\leq2^{\max\{q-1,0\}}(a^{q}+b^{q}).$$
\end{LemD}

\subsection{The proof of Theorem \ref{thm-2}}
We first prove $(1)\Rightarrow(2)$. By Lemma \ref{lem-1}, there
exist  functions
$f_{j}\in\mathscr{B}_{\omega}(\mathbb{D})\cap \mathscr{H}(\mathbb{D})$
$(j\in\{1,2\})$ such that, for $z\in\mathbb{D}$,
\be\label{eq-xp-1-1}
\sum_{j=1}^{2}|f_{j}'(z)|\geq\omega(|z|).\ee

Since
$C_{\phi}:~\mathscr{B}_{\omega}(\mathbb{D})\rightarrow\mathscr{PH}^{p}(\mathbb{B}^{n})$
is a bounded operator, by (\ref{g-1}),  we see that
\beq\label{eq-ch-0-1}\infty&>&\int_{\partial\mathbb{B}^{n}}\big(\mathscr{G}(C_{\phi}(f_{j}))(\zeta)\big)^{p}d\sigma(\zeta)\\
\nonumber &=&
\int_{\partial\mathbb{B}^{n}}\left(\int_{0}^{1}|f_{j}'(\phi(r\zeta))|^{2}|\nabla\,\phi(r\zeta)|^{2}(1-r)dr\right)^{\frac{p}{2}}d\sigma(\zeta)
\eeq
 for $j\in\{1,2\}$.
It follows from (\ref{eq-xp-1-1}), (\ref{eq-ch-0-1}) and Lemma D
that

\beqq
\infty&>&\sum_{j=1}^{2}\int_{\partial\mathbb{B}^{n}}\big(\mathscr{G}(C_{\phi}(f_{j}))(\zeta)\big)^{p}d\sigma(\zeta)\\
&\geq&\mathscr{M}_{p}\int_{\partial\mathbb{B}^{n}}\left(\int_{0}^{1}\left(\sum_{j=1}^{2}
|f_{j}'(\phi(r\zeta))|\right)^{2}|\nabla\,\phi(r\zeta)|^{2}(1-r)dr\right)^{\frac{p}{2}}d\sigma(\zeta)\\
&\geq&\mathscr{M}_{p}\int_{\partial\mathbb{B}^{n}}\left(\int_{0}^{1}
\omega^{2}(|\phi(r\zeta)|)|\nabla\,\phi(r\zeta)|^{2}(1-r)dr\right)^{\frac{p}{2}}d\sigma(\zeta),
\eeqq where $\mathscr{M}_{p}=2^{-p/2-\max\{p/2-1,0\}}$.

Next, we prove $(2)\Rightarrow(1)$. We split the proof of this case
into two steps.

$\mathbf{Step~1.}$ We first prove
$C_{\phi}(f)\in\mathscr{PH}^{p}(\mathbb{B}^{n})$. Since $\mathbb{D}$
is a simply connected domain, we see that $f$ admits the canonical
decomposition $f=f_{1}+\overline{f_{2}}$, where $f_{1}$ and $f_{2}$
are analytic in $\mathbb{D}$ with $f_{2}(0)=0$. For $j\in\{1,2\}$,
elementary calculations lead to \beq\label{rm-1}
\left|\nabla C_{\phi}(f_{j})(z)\right|&=&|f_{j}'(\phi(z))||\nabla\,\phi(z)|\leq\|f_{j}\|_{\mathscr{B}_{\omega}(\mathbb{D})}|\nabla\,\phi(z)|
\omega(|\phi(z)|)\\ \nonumber
&\leq&\|f\|_{\mathscr{B}_{\omega}(\mathbb{D})}|\nabla\,\phi(z)|\omega(|\phi(z)|).\eeq

For $j\in\{1,2\}$, let $$\mathscr{Y}_{j}=\int_{\partial\mathbb{B}^{n}}\big(\mathscr{G}(C_{\phi}(f_{j}))(\zeta)\big)^{p}d\sigma(\zeta).$$
Consequently,  there is a constant
$C=\|f\|_{\mathscr{B}_{\omega}(\mathbb{D})}^p$ such that

\beq \label{eq-ch-0h1}
\mathscr{Y}_{j}
&=&
\int_{\partial\mathbb{B}^{n}}\left(\int_{0}^{1}|f_{j}'(\phi(r\zeta))|^{2}|\nabla\,\phi(r\zeta)|^{2}(1-r)dr\right)^{\frac{p}{2}}d\sigma(\zeta)\\ \nonumber
&\leq&C\int_{\partial\mathbb{B}^{n}}\left(\int_{0}^{1}\omega^{2}(|\phi(r\zeta)|)|\nabla\,\phi(r\zeta)|^{2}(1-r)dr\right)^{\frac{p}{2}}d\sigma(\zeta)\\
\nonumber &<&\infty, \eeq which, together with (\ref{g-1}), implies
that $C_{\phi}(f_{j})\in\mathscr{H}^{p}(\mathbb{B}^{n})$ for
$j\in\{1,2\}$. Hence
$C_{\phi}(f)\in\mathscr{PH}^{p}(\mathbb{B}^{n})$.

$\mathbf{Step~2.}$ In this step, we will show that $C_{\phi}$ is a
bounded operator.  Without loss of generality, we assume that
$\|f\|_{\mathscr{B}_{\omega}(\mathbb{D})}\neq0$, and $f_{j}$ are not
constant functions for $j\in\{1,2\}$. Since, for $j\in\{1,2\}$,
\beq\label{eq-fgy-1}
|C_{\phi}(f_{j})(0)|&\leq&|f_{j}(0)|+|f_{j}(\phi(0))-f_{j}(0)|\\
\nonumber
&\leq&|f_{j}(0)|+|\phi(0)|\int_{0}^{1}|f_{j}'(\phi(0)t)|dt\\
\nonumber
&\leq&|f_{j}(0)|+\|f_{j}\|_{\mathscr{B}_{\omega}(\mathbb{D})}|\phi(0)|\int_{0}^{1}\omega(|\phi(0)|t)dt\\
\nonumber
&\leq&\|f_{j}\|_{\mathscr{B}_{\omega}(\mathbb{D})}\left(1+|\phi(0)|\omega(|\phi(0)|)\right),
\eeq by (\ref{g-1-1}) and (\ref{eq-ch-0h1}), we see that there is a
positive constant $C$ such that

\beq\label{c-h-1-3}
\|C_{\phi}(f_{j})\|_{p}^{p}&\leq&\,C\left(|C_{\phi}(f_{j})(0)|^{p}+\int_{\partial\mathbb{B}^{n}}\big(\mathscr{G}(C_{\phi}(f_{j}))(\zeta)\big)^{p}d\sigma(\zeta)\right)\\
\nonumber
&\leq&C\|f_{j}\|_{\mathscr{B}_{\omega}(\mathbb{D})}^{p}C(\phi),\eeq where

\beqq
C(\phi)&=&\big(1+|\phi(0)|\omega(|\phi(0)|)\big)^{p}\\
&&+\int_{\partial\mathbb{B}^{n}}
\left(\int_{0}^{1}\omega^{2}(|\phi(r\zeta)|)|\nabla\,\phi(r\zeta)|^{2}(1-r)dr\right)^{\frac{p}{2}}d\sigma(\zeta).\eeqq

Combining (\ref{c-h-1-3}) and Lemma D 
 gives

\beqq\|C_{\phi}(f)\|_{p}^{p}&\leq&2^{\max\{p-1,0\}}\left(\|C_{\phi}(f_{1})\|_{p}^{p}+\|C_{\phi}(f_{2})\|_{p}^{p}\right)\\
&\leq&2^{\max\{p-1,0\}}CC(\phi)\big(\|f_{1}\|_{\mathscr{B}_{\omega}(\mathbb{D})}^{p}+\|f_{2}\|_{\mathscr{B}_{\omega}(\mathbb{D})}^{p}\big)\\
&\leq&2^{1+\max\{p-1,0\}}CC(\phi)\|f\|^{p}_{\mathscr{B}_{\omega}(\mathbb{D})}.
\eeqq
Therefore,
$C_{\phi}:~\mathscr{B}_{\omega}(\mathbb{D})\rightarrow\mathscr{PH}^{p}(\mathbb{B}^{n})$
is a bounded operator.

Now, we prove $(1)\Rightarrow(3)$. By Lemma \ref{lem-1}, there
exist  functions
$f_{j}\in\mathscr{B}_{\omega}(\mathbb{D})\cap \mathscr{H}(\mathbb{D})$
$(j\in\{1,2\})$ such that, for $z\in\mathbb{D}$,
\be\label{eq-xp-1-11}\sum_{j=1}^{2}|f_{j}'(z)|\geq\omega(|z|).\ee

Since
$C_{\phi}:~\mathscr{B}_{\omega}(\mathbb{D})\rightarrow\mathscr{PH}^{p}(\mathbb{B}^{n})$
is a bounded operator, by (\ref{g-1-1.1}),  we see that
\beq\label{eq-ch-0-11}\nonumber 
\infty&>&\int_{\partial\mathbb{B}^{n}}\big(\mathscr{A}_{\alpha}C_{\phi}(f_{j})(\zeta)\big)^{p}d\sigma(\zeta)\\
&=&
\int_{\partial\mathbb{B}^{n}}\left(\int_{D_{\alpha}^{n}(\zeta)}|f_{j}'(\phi(z))|^{2}|\nabla\,\phi(z)|^{2}(1-|z|)^{1-n}
dV(z)\right)^{\frac{p}{2}}d\sigma(\zeta)
\eeq
 for $j\in\{1,2\}$. It follows from (\ref{eq-xp-1-11}), (\ref{eq-ch-0-11}) and Lemma D that

\beqq
\infty&>&\sum_{j=1}^{2}\int_{\partial\mathbb{B}^{n}}\big(\mathscr{A}_{\alpha}C_{\phi}(f_{j})(\zeta)\big)^{p}d\sigma(\zeta)\\
&\geq&\mathscr{M}_{p}\int_{\partial\mathbb{B}^{n}}\left(\int_{D_{\alpha}^{n}(\zeta)}\left(\sum_{j=1}^{2}
|f_{j}'(\phi(z)))|\right)^{2}\mathscr{V}(z)dV(z)\right)^{\frac{p}{2}}d\sigma(\zeta)\\
&\geq&\mathscr{M}_{p}\int_{\partial\mathbb{B}^{n}}\left(\int_{D_{\alpha}^{n}(\zeta)}
\omega^{2}(|\phi(z)|)\mathscr{V}(z)dV(z)\right)^{\frac{p}{2}}d\sigma(\zeta),
\eeqq  where $\mathscr{V}(z)=|\nabla\,\phi(z)|^{2}(1-|z|)^{1-n}$.


The proof of  $(3)\Rightarrow(1)$ is similar to the proof of  $(2)\Rightarrow(1)$. We only need to replace ``(\ref{g-1-1})" and ``$\int_{\partial\mathbb{B}^{n}}\big(\mathscr{G}(C_{\phi}(f_{j}))(\zeta)\big)^{p}d\sigma(\zeta)$" by ``(\ref{g-1-1.2})" and $$``\int_{\partial\mathbb{B}^{n}}\big(\mathscr{A}_{\alpha}C_{\phi}(f_{j})(\zeta)\big)^{p}d\sigma(\zeta)",$$ respectively, in the proof of $(2)\Rightarrow(1)$,
where $j\in\{1,2\}$.

At last, we prove $(1)\Leftrightarrow(4)$. Since $(4)\Rightarrow(1)$
is obvious, we only need to prove $(1)\Rightarrow(4)$. Let
$C_{\phi}:~\mathscr{B}_{\omega}(\mathbb{D})\rightarrow\mathscr{PH}^{p}(\mathbb{B}^{n})$
be a bounded operator. Then
$f\circ\phi\in\mathscr{PH}^{p}(\mathbb{B}^{n})$ for all
$f\in\mathscr{B}_{\omega}(\mathbb{D})$, and
$$|\phi(\zeta)|:=|\lim_{r\rightarrow1^{-}}\phi(r\zeta)|\leq1$$ exists for almost every $\zeta\in\partial\mathbb{B}^{n}$.
Suppose that $\{f_{k}=h_{k}+\overline{g_{k}}\}$ is a sequence of
$\mathscr{B}_{\omega}(\mathbb{D})$ such that
$\|f_{k}\|_{\mathscr{B}_{\omega}(\mathbb{D})}\leq1$, where $h_{k}$ and
$g_{k}$ are analytic in $\mathbb{D}$ with $g_{k}(0)=0$. Next, we
will show that $\{C_{\phi}(f_{k})\}$ has a convergent subsequence in
$\mathscr{PH}^{p}(\mathbb{B}^{n})$. Since \beqq
\max\{ |h_{k}(z)|, |g_k(z)|\}
&\leq&|f_{k}(0)|+\|f_{k}\|_{\mathscr{B}_{\omega}(\mathbb{D})}\int_{0}^{1}\omega(t|z|)|z|dt\\
&\leq&1+\int_{0}^{1}\omega(t|z|)|z|dt
\\
&\leq &
1+\omega(|z|)
<\infty,
\eeqq
$\{ h_{k}\}$ and $\{ g_k\}$ are normal families.
Then there
are subsequences of $\{ h_{k}\}$ and $\{ g_k\}$ that converge uniformly on compact
subsets of $\mathbb{D}$ to holomorphic functions $h$ ang $g$, respectively.
Without loss of generality, we assume that the
sequences $\{ h_{k}\}$ and $\{ g_{k}\}$ themselves converge to $h$ and $g$, respectively.
Then the
sequence $f_{k}=h_{k}+\overline{g_{k}}$ itself converges uniformly on compact
subsets of $\mathbb{D}$
to the harmonic function
$f=h+\overline{g}$ with $g(0)=0$.
Consequently,
$$1\geq\lim_{k \rightarrow\infty}\left\{|f_{k}(0)|+\mathscr{B}_{\omega}^{f_{k}}(z)\right\}=|f(0)|+\mathscr{B}_{\omega}^{f}(z),$$
which gives that $f\in\mathscr{B}_{\omega}(\mathbb{D})$ with
$\|f\|_{\mathscr{B}_{\omega}(\mathbb{D})}\leq 1$. It follows from
$(1)\Rightarrow(2)$ that
\begin{equation}
\label{dominate}
\int_{0}^{1}\omega^{2}(|\phi(r\zeta)|)|\nabla \phi(r\zeta)|^{2}(1-r)dr<\infty,
\quad
\mbox{a.e. }
\zeta\in \partial \mathbb{B}^n
\end{equation}
and
\begin{equation}
\label{dominate2}
\int_{ \partial \mathbb{B}^n}\left(\int_{0}^{1}|\nabla \phi(r\zeta)|^{2}\omega^{2}(|\phi(r\zeta)|)(1-r)
dr\right)^{\frac{p}{2}}\,d\sigma(\zeta)<\infty.
\end{equation}
Also, we have
\beq\label{eq-x-1}
\lefteqn{\left|(h_{k}-h)'\circ\phi(r\zeta)\right|^{2}|\nabla\phi(r\zeta)|^{2}(1-r)}
\\
&&\quad \leq
\|h_{k}-h\|_{\mathscr{B}_{\omega}(\mathbb{D})}^{2}
\omega^{2}(|\phi(r\zeta)|)|\nabla \phi(r\zeta)|^{2}(1-r)
\nonumber
\\
&&\quad \leq 4\omega^{2}(|\phi(r\zeta)|)|\nabla\phi(r\zeta)|^{2}(1-r)
\nonumber
\eeq
for $\zeta\in  \partial \mathbb{B}^n$.

Applying (\ref{dominate}), (\ref{dominate2}), (\ref{eq-x-1})  and the control convergence theorem twice,
we have

\beqq &&\lim_{k\rightarrow\infty}\int_{\partial\mathbb{B}^{n}}
\left(\int_{0}^{1}\left|(h_{k}-h)'\circ\phi(r\zeta)\right|^{2}|\nabla\,\phi(r\zeta)|^{2}(1-r)dr\right)^{\frac{p}{2}}d\sigma(\zeta)\\
&=&\int_{\partial\mathbb{B}^{n}}
\left(\int_{0}^{1}\lim_{k\rightarrow\infty}\left|(h_{k}-h)'\circ\phi(r\zeta)\right|^{2}|\nabla\,\phi(r\zeta)|^{2}(1-r)dr\right)^{\frac{p}{2}}d\sigma(\zeta)=0,
\eeqq which, together with (\ref{g-1-1}), implies that

$$\lim_{k\rightarrow\infty}\|h_{k}\circ\phi-h\circ \phi\|_p=0. $$
Consequently,
 \be\label{eq-x-2}C_{\phi}(h_{k})\rightarrow\,C_{\phi}(h)
 \quad \mbox{in \ }\mathscr{PH}^{p}(\mathbb{B}^{n})\ee
 as $k\rightarrow\infty$.
By using similar reasoning as in the proof of (\ref{eq-x-2}), we
have
\be\label{eq-x-3}C_{\phi}(g_{k})\rightarrow\,C_{\phi}(g)\quad \mbox{in \ }\mathscr{PH}^{p}(\mathbb{B}^{n})\ee
 as $k\rightarrow\infty$.
Therefore, by (\ref{eq-x-2}) and (\ref{eq-x-3}), we conclude that
$C_{\phi}(f)=C_{\phi}(h)+\overline{C_{\phi}(g)}\in\mathscr{PH}^{p}(\mathbb{B}^{n})$,
and $C_{\phi}(f_{k})\rightarrow\,C_{\phi}(f)$ in $\mathscr{PH}^{p}(\mathbb{B}^{n})$
as $k\rightarrow\infty$. The proof of this theorem is finished.\qed

\begin{ThmE}{\rm(\cite[Theorem 1.1]{Pav})}\label{Thm-P-1}
Let $p\in(0,\infty)$ and $f$ be  a holomorphic function in
$\mathbb{D}$. Then the followings are equivalent:
\begin{enumerate}
\item[{\rm (1)}] $f\in\mathscr{H}^{p}(\mathbb{D})$;

\item[{\rm (2)}] $\mathcal{G}[f]\in L^p(\mathbb{T}),$
where
\[
\mathcal{G}[f](\zeta)=\left(\int_{0}^{1}(1-r)|f'(r\zeta)|^{2}dr\right)^{\frac{1}{2}},
\quad
\zeta \in \mathbb{T};
\]

\item[{\rm (3)}]
$\mathcal{G}_*[f]\in L^p(\mathbb{T}),$
where
\[
\mathcal{G}_*[f](\zeta)=\left(\int_{0}^{1}(1-r)\sup_{\rho\in(0,r)}|f'(\rho\,\zeta)|^{2}dr\right)^{\frac{1}{2}},
\quad
\zeta \in \mathbb{T};
\]

\item[{\rm (4)}]
$\mathcal{G}_d[f]\in L^p(\mathbb{T}),$
where
\[
\mathcal{G}_d[f](\zeta)=\left(\sum_{k=0}^{\infty}2^{-2k}|f'(r_{k}\zeta)|^{2}\right)^{\frac{1}{2}},
\quad
\zeta \in \mathbb{T}
\]
and $r_{k}=1-2^{-k}$.
\end{enumerate}
Furthermore, there are constants
$C_1$, $C_2$, $C_3$ and $C_4$ independent of $f$ such that
\[
\| f-f(0)\|_p\leq C_1\|\mathcal{G}[f]\|_p
\leq
C_2\|\mathcal{G}_*[f]\|_p
\leq
C_3\|\mathcal{G}_d[f]\|_p
\leq
C_4\| f-f(0)\|_p.
\]
\end{ThmE}

\subsection{The proof of Theorem \ref{thm-3}} Since
$(1)\Leftrightarrow(2)\Leftrightarrow(5)\Leftrightarrow(6)$ easily follows from Theorem
\ref{thm-2}, we only need to prove
$(1)\Leftrightarrow(3)\Leftrightarrow(4)$. We first prove
$(1)\Rightarrow(3)$. By Lemma \ref{lem-1}, there are two analytic
functions
$f_{j}\in\mathscr{B}_{\omega}(\mathbb{D})\cap \mathscr{H}(\mathbb{D})$
$(j\in\{1,2\})$ such that, for $z\in\mathbb{D}$,
\be\label{eq-xpj-1-1}\sum_{j=1}^{2}|f_{j}'(z)|\geq\omega(|z|).\ee

Since
$C_{\phi}:~\mathscr{B}_{\omega}(\mathbb{D})\rightarrow\mathscr{PH}^{p}(\mathbb{D})$
is a bounded operator, by Theorem E,
we see that
\beq\label{eq-ch-1-7}\infty&>&\frac{1}{2\pi}\int_{0}^{2\pi}\left(\sum_{k=0}^{\infty}2^{-2k}\left|\big(f_{j}(\phi(r_{k}e^{i\theta})\big)'\right|^{2}\right)^{\frac{p}{2}}d\theta\\
\nonumber &=&
\frac{1}{2\pi}\int_{0}^{2\pi}\left(\sum_{k=0}^{\infty}2^{-2k}\left|f_{j}'(\phi(r_{k}e^{i\theta}))\phi'(r_{k}e^{i\theta})\right|^{2}\right)^{\frac{p}{2}}d\theta
\eeq
 for $j\in\{1,2\}$.
It follows from (\ref{eq-xpj-1-1}), (\ref{eq-ch-1-7}) and Lemma D
 that

\beqq
\infty&>&\sum_{j=1}^{2}\int_{0}^{2\pi}\left(\sum_{k=0}^{\infty}2^{-2k}
\left|\big(f_{j}(\phi(r_{k}e^{i\theta})\big)'\right|^{2}\right)^{\frac{p}{2}}\frac{d\theta}{2\pi}\\
&\geq&\mathscr{M}_{p}\int_{0}^{2\pi}\left(\sum_{k=0}^{\infty}2^{-2k}|\phi'(r_{k}e^{i\theta})|^{2}\left(\sum_{j=1}^{2}|f_{j}'(\phi(r_{k}e^{i\theta})|\right)^{2}\right)^{\frac{p}{2}}\frac{d\theta}{2\pi}\\
&\geq&\mathscr{M}_{p}\int_{0}^{2\pi}\left(\sum_{k=0}^{\infty}2^{-2k}|\phi'(r_{k}e^{i\theta})|^{2}\omega^{2}(|\phi(r_{k}e^{i\theta})|)\right)^{\frac{p}{2}}\frac{d\theta}{2\pi},
\eeqq where $\mathscr{M}_{p}$ is defined in the proof of Theorem
\ref{thm-2}.

Next, we prove $(3)\Rightarrow(1)$. We first show
$C_{\phi}(f)\in\mathscr{PH}^{p}(\mathbb{D})$. Let
$f\in\mathscr{B}_{\omega}(\mathbb{D})$. Since $\mathbb{D}$ is a simply
connected domain, we see that $f$ admits the canonical decomposition
$f=f_{1}+\overline{f_{2}}$, where $f_{1}$ and $f_{2}$ are analytic
in $\mathbb{D}$ with $f_{2}(0)=0$. For $j\in\{1,2\}$,
let 
$$\mathscr{Z}_{j}=\int_{0}^{2\pi}\left(\sum_{k=0}^{\infty}\frac{\left|\big(f_{j}(\phi(r_{k}e^{i\theta})\big)'\right|^{2}}{2^{2k}}\right)^{\frac{p}{2}}\frac{d\theta}{2\pi}.$$
Then,
 by
(\ref{rm-1}), we see that there is a constant
$C=\|f\|_{\mathscr{B}_{\omega}(\mathbb{D})}^p$ such that

\beqq
\mathscr{Z}_{j}
&=&
\int_{0}^{2\pi}\left(\sum_{k=0}^{\infty}\frac{\left|f_{j}'(\phi(r_{k}e^{i\theta}))\phi'(r_{k}e^{i\theta})\right|^{2}}{2^{2k}}\right)^{\frac{p}{2}}\frac{d\theta}{2\pi}\\
&\leq&C
\int_{0}^{2\pi}\left(\sum_{k=0}^{\infty}\frac{\left|\omega(|\phi(r_{k}e^{i\theta})|)\phi'(r_{k}e^{i\theta})\right|^{2}}{2^{2k}}\right)^{\frac{p}{2}}\frac{d\theta}{2\pi}\\
&<&\infty, \eeqq which, together with Theorem E, 
implies
that $C_{\phi}(f_{j})\in\mathscr{PH}^{p}(\mathbb{D})$. Hence
$$C_{\phi}(f)=C_{\phi}(f_{1})+\overline{C_{\phi}(f_{2})}\in\mathscr{PH}^{p}(\mathbb{D}).$$

Now we come to prove $C_{\phi}$ is a bounded operator. By
(\ref{eq-fgy-1}) and Theorem E,
we see that there is a
positive constant $C$, depending only on $p$, such that, for
$j\in\{1,2\}$,

\beq\label{eq-fgy-3}
\nonumber
\|C_{\phi}(f_{j})\|_{p}^{p}&\leq&\,C
\left(|f_j(\phi(0))|^p+
\int_{0}^{2\pi}\left(
\sum_{k=0}^{\infty}\frac{\left|
\big(f_{j}(\phi(r_{k}e^{i\theta})\big)'\right|^{2}}{2^{2k}}\right)^{\frac{p}{2}}\frac{d\theta}{2\pi}\right)\\
&\leq&C\|f_{j}\|_{\mathscr{B}_{\omega}(\mathbb{D})}^pC^{\ast}(\phi), \eeq
where $$C^{\ast}(\phi)=\left(1+|\phi(0)|\omega(|\phi(0)|)\right)^p+\int_{0}^{2\pi}\left(\sum_{k=0}^{\infty}\frac{\left|\omega(|\phi(r_{k}e^{i\theta})|)
\phi'(r_{k}e^{i\theta})\right|^{2}}{2^{2k}}\right)^{\frac{p}{2}}\frac{d\theta}{2\pi}.$$
It follows from (\ref{eq-fgy-3}) and Lemma D  
that

\beqq\|C_{\phi}(f)\|_{p}^{p}&\leq&2^{\max\{p-1,0\}}\left(\|C_{\phi}(f_{1})\|_{p}^{p}+\|C_{\phi}(f_{2})\|_{p}^{p}\right)\\
&\leq&2^{\max\{p-1,0\}}CC^{\ast}(\phi)\big(\|f_{1}\|_{\mathscr{B}_{\omega}(\mathbb{D})}^{p}+\|f_{2}\|_{\mathscr{B}_{\omega}(\mathbb{D})}^{p}\big)\\
&\leq&2^{1+\max\{p-1,0\}}CC^{\ast}(\phi)\|f\|^{p}_{\mathscr{B}_{\omega}(\mathbb{D})},
\eeqq which implies that $C_{\phi}$ is a bounded operator.

Now we prove $(1)\Rightarrow(4)$. By Lemma \ref{lem-1}, there are
two analytic functions
$f_{j}\in\mathscr{B}_{\omega}(\mathbb{D})\cap \mathscr{H}(\mathbb{D})$
$(j\in\{1,2\})$ such that, for $z\in\mathbb{D}$,
\be\label{eq-xpj-1-5}\sum_{j=1}^{2}|f_{j}'(z)|\geq\omega(|z|).\ee
Since $C_{\phi}$ is a bounded operator, by (\ref{eq-xpj-1-5}),
Theorem E 
and Lemma D, 
we see that

\beqq
\infty&>&\frac{1}{\mathscr{M}_{p}}\sum_{j=1}^{2}\int_{0}^{2\pi}\left(\int_{0}^{1}\delta(r)\sup_{\rho\in(0,r)}
\left|\big(f_{j}(\phi(\rho\,e^{i\theta})\big)'\right|^{2}dr\right)^{\frac{p}{2}}\frac{d\theta}{2\pi}\\
&\geq&\int_{0}^{2\pi}\left(\int_{0}^{1}\delta(r)\sup_{\rho\in(0,r)}\left(
|\phi'(\rho\,e^{i\theta})|^{2}\left(\sum_{j=1}^{2}|f_{j}'(\phi(\rho e^{i\theta})|\right)^{2}\right)dr\right)^{\frac{p}{2}}\frac{d\theta}{2\pi}\\
&\geq&\int_{0}^{2\pi}\left(\int_{0}^{1}\delta(r)\sup_{\rho\in(0,r)}\left(
|\phi'(\rho\,e^{i\theta})|^{2}\omega^{2}(|\phi(\rho\,e^{i\theta})|)\right)dr\right)^{\frac{p}{2}}\frac{d\theta}{2\pi},
\eeqq where $\delta(r)=1-r$ and $\mathscr{M}_{p}$ is defined in the proof of Theorem
\ref{thm-2}.

At last, we prove $(4)\Rightarrow(1)$. Let
$f\in\mathscr{B}_{\omega}(\mathbb{D})$. It is not difficult to know that,
for $z\in\mathbb{D}$,
  $f$  has the canonical decomposition
$f(z)=f_{1}(z)+\overline{f_{2}(z)}$ , where $f_{1}$ and $f_{2}$ are
analytic in $\mathbb{D}$ with $f_{2}(0)=0$. Then, by (\ref{rm-1}),
we see that

\beq\label{qew-1}\nonumber
\infty&>&M_{j}\int_{0}^{2\pi}\left(\int_{0}^{1}\delta(r)\sup_{\rho\in(0,r)}\left(
|\phi'(\rho\,e^{i\theta})|^{2}\omega^{2}(|\phi(\rho\,e^{i\theta})|)\right)dr\right)^{\frac{p}{2}}\frac{d\theta}{2\pi}\\
\nonumber
&\geq&\int_{0}^{2\pi}\left(\int_{0}^{1}\delta(r)\sup_{\rho\in(0,r)}\left(
|\phi'(\rho\,e^{i\theta})|^{2}|f_{j}'(\phi(\rho e^{i\theta})|^{2}\right)dr\right)^{\frac{p}{2}}\frac{d\theta}{2\pi}\\ 
&=&\int_{0}^{2\pi}\left(\int_{0}^{1}\delta(r)\sup_{\rho\in(0,r)}
\left|\big(f_{j}(\phi(\rho\,e^{i\theta})\big)'\right|^{2}dr\right)^{\frac{p}{2}}\frac{d\theta}{2\pi},
 \eeq
where  $j\in\{1,2\}$ and $M_{j}=\|f_{j}\|_{\mathscr{B}_{\omega}(\mathbb{D})}^p$. It follows from (\ref{qew-1}) and Theorem E
 that
$C_{\phi}(f_{j})\in\mathscr{PH}^{p}(\mathbb{D})~(j=1,2)$, which
implies that
$$C_{\phi}(f)=C_{\phi}(f_{1})+\overline{C_{\phi}(f_{2})}\in\mathscr{PH}^{p}(\mathbb{D}).$$

Next, we prove $C_{\phi}$ is also a bounded operator. By
(\ref{eq-fgy-1}) and Theorem  E, 
we see that there is a
positive constant $C$, depending only on $p$, such that, for
$j\in\{1,2\}$,

\beq\label{eq-fgy-6}
\|C_{\phi}(f_{j})\|_{p}^{p}
&\leq&\,C
\Bigg( |f_j(\phi(0))|^p\\ \nonumber
&&+
\int_{0}^{2\pi}\left(\int_{0}^{1}\delta(r)\sup_{\rho\in(0,r)}
\left|\big(f_{j}(\phi(\rho\,e^{i\theta})\big)'\right|^{2}dr\right)^{\frac{p}{2}}\frac{d\theta}{2\pi}\Bigg)\\
\nonumber &\leq&CMM_{j}, \eeq where
\begin{eqnarray*}
\lefteqn{M=\left(1+|\phi(0)|\omega(|\phi(0)|)\right)^p}
\\
&& \quad
+\int_{0}^{2\pi}\left(\int_{0}^{1}(1-r)\sup_{\rho\in(0,r)}\left(
|\phi'(\rho\,e^{i\theta})|\omega(|\phi(\rho\,e^{i\theta})|)\right)^{2}dr\right)^{\frac{p}{2}}\frac{d\theta}{2\pi}.
\end{eqnarray*}
It follows from (\ref{eq-fgy-6}) and Lemma D 
that

\beqq\|C_{\phi}(f)\|_{p}^{p}&\leq&2^{\max\{p-1,0\}}\left(\|C_{\phi}(f_{1})\|_{p}^{p}+\|C_{\phi}(f_{2})\|_{p}^{p}\right)\\
&\leq&2^{\max\{p-1,0\}}CM\big(M_{1}+M_{2}\big)\\
&\leq&2^{1+\max\{p-1,0\}}\|f\|^{p}_{\mathscr{B}_{\omega}(\mathbb{D})}CM,
\eeqq which implies that
$C_{\phi}:~\mathscr{B}_{\omega}(\mathbb{D})\rightarrow\mathscr{PH}^{p}(\mathbb{D})$
is a bounded operator. The proof of this theorem is finished.
 \qed

\subsection{The proof of Proposition \ref{p-1}} If $k=1$, then we can take $r_1=0$.
Since, for $k\in\{2,3,\ldots\}$,
$$\lim_{x\rightarrow0^{+}}\mu_{\omega,k}(x)=\lim_{x\rightarrow1^{-}}\mu_{\omega,k}(x)=0,$$
 we see that $\mathscr{E}_{\omega}(k)$ is not an empty set. 
Consequently, for $k\in\{2,3,\ldots\}$, we can choose  points $r_{k}>0$ and $r_{k+1}>0$ from the sets $\mathscr{E}_{\omega}(k)$ and $\mathscr{E}_{\omega}(k+1)$, respectively,
such that $\mu_{\omega,k}(r_{k})\geq \mu_{\omega,k}(r_{k+1})$ which  is equivalent to
\be\label{w-1}\frac{\omega(r_{k+1})}{\omega(r_{k})}\geq \frac{r_{k+1}^{k-1}}{r_{k}^{k-1}}.\ee
On the other hand, for $k\in\{2,3,\ldots\}$, we have $\mu_{\omega,k+1}(r_{k+1})\geq\mu_{\omega,k+1}(r_{k})$ which implies that
\be\label{w-2}\frac{\omega(r_{k+1})}{\omega(r_{k})}\leq\frac{r_{k+1}^{k}}{r_{k}^{k}}.\ee
Combining (\ref{w-1}) and (\ref{w-2}) gives $r_{k}\leq \,r_{k+1}$.
\qed

\subsection{The proof of Theorem \ref{Thm-Bloch}}
(1)
We first prove the sufficiency. By (\ref{eqx-4}), we have
\be\label{eqx-5}\|\phi\|_{\mathscr{B}_{\omega_{2}}(\mathbb{B}^n)}<\infty.\ee
For  $f\in\mathscr{B}_{\omega_{1}}(\mathbb{D})$,
we have
\begin{eqnarray*}
|C_{\phi}(f)(0)|&\leq&|f(0)|+|f(\phi(0))-f(0)|\\
\nonumber
&\leq&|f(0)|+|\phi(0)|\int_{0}^{1}\Lambda_f(\phi(0)t)dt\\
\nonumber
&\leq&|f(0)|+\|f\|_{\mathscr{B}_{\omega_1}(\mathbb{D})}|\phi(0)|\int_{0}^{1}\omega_1(|\phi(0)|t)dt\\
\nonumber
&\leq&\|f\|_{\mathscr{B}_{\omega_1}(\mathbb{D})}\left(1+|\phi(0)|\omega_1(|\phi(0)|)\right).
\end{eqnarray*}
So, it suffices to show that there exists a constant independent of $f$ such that
\[
\| f\circ \phi\|_{\mathscr{B}_{\omega_2}(\mathbb{B}^n),s}\leq C \|f\|_{\mathscr{B}_{\omega_1}(\mathbb{D}),s},
\quad
f\in\mathscr{B}_{\omega_{1}}(\mathbb{D}).
\]
We split the remaining proof into two cases.

\noindent $\mathbf{Case~1.}$ If $$\sup_{z\in\mathbb{B}^n}|\phi(z)|<1,$$ then there is a constant $\rho_{0}\in(0,1)$ such that
\be\label{eqx-6}\sup_{z\in\mathbb{B}^n}|\phi(z)|<\rho_{0}.\ee
For  $f\in\mathscr{B}_{\omega_{1}}(\mathbb{D})$, it follows from (\ref{eqx-5}) and (\ref{eqx-6}) that

\beqq
\sup_{z\in\mathbb{B}^n}\left\{\frac{\Lambda_f(\phi(z))\|D\phi(z)\|}{\omega_{2}(|z|)}\right\}&=&
\sup_{z\in\mathbb{B}^n}\left\{\frac{\Lambda_f(\phi(z))}{\omega_{1}(|\phi(z)|)}\omega_{1}(|\phi(z)|)\frac{\|D\phi(z)\|}{\omega_{2}(|z|)}\right\}\\
&\leq&\|f\|_{\mathscr{B}_{\omega_1}(\mathbb{D}),s}\omega_{1}(\rho_0)\sup_{z\in\mathbb{B}^n}
\left\{\frac{\|D\phi(z)\|}{\omega_{2}(|z|)}\right\}\\
&<&\infty,
\eeqq
which implies that $C_{\phi}:~\mathscr{B}_{\omega_{1}}(\mathbb{D})\rightarrow\mathscr{B}_{\omega_{2}}(\mathbb{B}^n)$ is bounded.

\noindent $\mathbf{Case~2.}$ If \be\label{eqx-7}\sup_{z\in\mathbb{B}^n}|\phi(z)|=1,\ee then, for $k\in\{1,2,\ldots\}$, let
$$\Omega_{k}=\left\{z\in\mathbb{B}^n:~r_{k}\leq|\phi(z)|\leq\, r_{k+1}\right\},$$ where $r_{1}=0$.
Since $\{r_k\}$ is a non-decreasing sequence satisfying \be\label{eqx-8.0} \lim_{k\rightarrow\infty}r_{k}^{k-1}=\gamma>0,\ee
we see that \be\label{eqx-8}\lim_{k\rightarrow\infty}r_{k}=1.\ee
Let $m$ be the smallest positive integer such that $\Omega_{m}$ is not an empty set. 
By (\ref{eqx-7}) and (\ref{eqx-8}), we see that $\Omega_{k}$ is not empty for every integer $k\in\{m,m+1,\ldots\}$, and
$\mathbb{B}^n=\cup_{k=m}^{\infty}\Omega_{k}.$
Then, for $k\in\{m,m+1,\ldots\}$, we have

\be\label{eqx-9}
\min_{z\in\Omega_{k}}\left\{\frac{\omega_{1}(r_{k})|\phi(z)|^{k-1}}{\omega_{1}(|\phi(z)|)}\right\}\geq
\left\{\frac{\omega_{1}(r_{k})r_{k}^{k-1}}{\omega_{1}(r_{k+1})}\right\}.
\ee
It follows from
 $\mu_{\omega_{1},k}(r_{k})\geq \mu_{\omega_{1},k}(r_{k+1})$  and  $\mu_{\omega_{1},k+1}(r_{k+1})\geq\mu_{\omega_{},k+1}(r_{k})$ that
\begin{equation}
\label{omega-estimate}
\frac{r_{k}^{k}}{r_{k+1}^{k}}\leq\frac{\omega_{1}(r_{k})}{\omega_{1}(r_{k+1})}\leq\frac{r_{k}^{k-1}}{r_{k+1}^{k-1}},
\end{equation}
which, together with (\ref{eqx-8.0}) and (\ref{eqx-8}), yields that
\be\label{eqx-10}\lim_{k\rightarrow\infty}\frac{\omega_{1}(r_{k})}{\omega_{1}(r_{k+1})}=1.\ee
Combining (\ref{eqx-9}) and (\ref{eqx-10}) gives

\beqq
\lim_{k\rightarrow\infty}\min_{z\in\Omega_{k}}\left\{\frac{\omega_{1}(r_{k})|\phi(z)|^{k-1}}{\omega_{1}(|\phi(z)|)}\right\}\geq
\lim_{k\rightarrow\infty}\left\{\frac{\omega_{1}(r_{k})r_{k}^{k-1}}{\omega_{1}(r_{k+1})}\right\}=\gamma.
\eeqq
Hence there is a positive integer $m_{0}\geq\,m$ such that, for all $k\in\{m_{0},m_{0}+1,\ldots\}$,
\be\label{rrg}\min_{z\in\Omega_{k}}\left\{\frac{\omega_{1}(r_{k})|\phi(z)|^{k-1}}{\omega_{1}(|\phi(z)|)}\right\}\geq\frac{\gamma}{2},\ee
which implies that, for $f\in\mathscr{B}_{\omega_{1}}(\mathbb{D})$,
\beq\label{ef-1-1}
\|C_{\phi}(f)\|_{\mathscr{B}_{\omega_2}(\mathbb{B}^n),s}&\leq &\sup_{z\in\mathbb{B}^n}\left\{\frac{\Lambda_f(\phi(z))\|D\phi(z)\|}{\omega_{2}(|z|)}\right\}
\\ \nonumber
&=&
\sup_{k\geq\,m}\sup_{z\in\Omega_{k}}\left\{\frac{\Lambda_f(\phi(z))\|D\phi(z)\|}{\omega_{2}(|z|)}\right\}\\ \nonumber
&=&
\max\{ J_1(f,m_0), J_2(f,m_0)\},
\eeq
where
\[
J_1(f,m_0)=
\sup_{k\geq\,m_0}\sup_{z\in\Omega_{k}}\left\{\frac{\Lambda_f(\phi(z))\|D\phi(z)\|}{\omega_{2}(|z|)}\right\}
\]
and
\[
J_2(f,m_0)=
\sup_{m_0-1\geq k\geq\,m}\sup_{z\in\Omega_{k}}\left\{\frac{\Lambda_f(\phi(z))\|D\phi(z)\|}{\omega_{2}(|z|)}\right\}.
\]
By (\ref{rrg}), we see that
\beqq
J_1(f,m_0)
&=&
\sup_{k\geq\,m_0}\sup_{z\in\Omega_{k}}\frac{\frac{\Lambda_f(\phi(z))\|D\phi(z)\|\omega_{1}(r_{k})|\phi(z)|^{k-1}}{\omega_{1}(|\phi(z)|)}}
{\frac{\omega_{2}(|z|)\omega_{1}(r_{k})|\phi(z)|^{k-1}}{\omega_{1}(|\phi(z)|)}}\\
&\leq&\frac{2}{\gamma}\sup_{k\geq\,m_0}\sup_{z\in\Omega_{k}}\frac{\omega_{1}(r_{k})}{k}
\frac{\|D(\phi^{k}(z))\|\Lambda_f(\phi(z))}{\omega_{1}(|\phi(z)|)\omega_{2}(|z|)}\\
&\leq&\frac{2}{\gamma}\sup_{k\geq1}\frac{\omega_{1}(r_{k})}{k}\|\phi^{k}\|_{\mathscr{B}_{\omega_{2}}(\mathbb{B}^n)}\|f\|_{\mathscr{B}_{\omega_1}(\mathbb{D}),s}\\
&<&\infty
\eeqq
and
\[
J_2(f,m_0)\leq \omega_1(r_{m_0})\|\phi\|_{\mathscr{B}_{\omega_{2}}(\mathbb{B}^n)}\|f\|_{\mathscr{B}_{\omega_1}(\mathbb{D}),s},
\]
which imply that $C_{\phi}:~\mathscr{B}_{\omega_{1}}(\mathbb{D})\rightarrow\mathscr{B}_{\omega_{2}}(\mathbb{B}^n)$ is bounded.

Next, we begin to prove the necessity. For $k\in\{2,3,\ldots\}$, let $f(w)=w^{k},~w\in\mathbb{D}$.
Since $$\|f\|_{\mathscr{B}_{\omega_{1}}(\mathbb{D})}=\sup_{w\in\mathbb{D}}
\left\{\frac{k|w|^{k-1}}{\omega_{1}(|w|)}\right\}=\frac{kr_{k}^{k-1}}{\omega_{1}(r_{k})},$$
we see that

 \be\label{eqx-12}\lim_{k\rightarrow\infty}\frac{\|f\|_{\mathscr{B}_{\omega_{1}}(\mathbb{D})}\omega_{1}(r_{k})}{k}
=\lim_{k\rightarrow\infty}r_{k}^{k-1}=\gamma>0.\ee
From (\ref{eqx-12}), we see that there is a positive constant, independent of $k$, such that
$$\|f\|_{\mathscr{B}_{\omega_{1}}(\mathbb{D})}\leq\,C\frac{k}{\omega_{1}(r_{k})}$$
which is equivalent to $$\frac{1}{\|f\|_{\mathscr{B}_{\omega_{1}}(\mathbb{D})}}\geq\frac{\omega_{1}(r_{k})}{kC}.$$
For $w\in\mathbb{D}$, let $F(w)=w^{k}/\|w^{k}\|_{\mathscr{B}_{\omega_{1}}(\mathbb{D})}$. Then $\|F\|_{\mathscr{B}_{\omega_{1}}(\mathbb{D})}=1$.
Therefore, $$\infty>\|C_{\phi}\|\geq\|C_{\phi}(F)\|_{\mathscr{B}_{\omega_{2}}(\mathbb{B}^n)}=
\frac{\|\phi^{k}\|_{\mathscr{B}_{\omega_{2}}(\mathbb{B}^n)}}{\|f\|_{\mathscr{B}_{\omega_{1}}(\mathbb{D})}}
\geq\frac{\omega_{1}(r_{k})\|\phi^{k}\|_{\mathscr{B}_{\omega_{2}}(\mathbb{B}^n)}}{kC},$$
which implies (\ref{eqx-4}) is true.

(2)
Assume that $C_{\phi}$ is compact.
For $k\in\{1,2,\ldots\}$ and $w\in\mathbb{D}$, let $F_k(w)=w^{k}/\|w^{k}\|_{\mathscr{B}_{\omega_{1}}(\mathbb{D})}$. Then $\|F_k\|_{\mathscr{B}_{\omega_{1}}(\mathbb{D})}=1$.
For $k\in\{2,3,\ldots\}$ and  $m\in\{1,2,\ldots\}$, it follows from  (\ref{omega-estimate}) that
\[
\frac{\omega_1(r_k)}{\omega_1(r_{k+m})}\geq
\frac{r_k^kr_{k+1}\cdots r_{k+m-1}}{r_{k+m}^{k+m-1}},
\]
which implies that
\begin{equation}
\label{omega-estimate2}
\frac{\omega_1(r_{k+m})r_k^{k+m}}{k+m}
\leq
\frac{\omega_1(r_k)r_k^mr_{k+m}^{k+m-1}}{(k+m)r_{k+1}\cdots r_{k+m-1}}
\leq
\frac{\omega_1(r_k)r_kr_{k+m}^{k+m-1}}{k+m}.
\end{equation}
For $k\in\{2,3,\ldots\}$ and $w\in \mathbb{D}$, elementary calculations lead to
$$
F_k(w)=\frac{\omega_1(r_k)w^k}{kr_k^{k-1}},
$$
which, together with (\ref{eqx-8.0}),  (\ref{eqx-8}) and (\ref{omega-estimate2}),
yields that
$F_k\to 0$ locally uniformly in $\mathbb{D}$ as $k\to\infty$.
Since $C_{\phi}$ is compact,
we deduce that
\[
\lim_{k\to \infty}\frac{\omega_1(r_k)}{kr_k^{k-1}}\|\phi^{k}\|_{\mathscr{B}_{\omega_{2}}(\mathbb{B}^n)}
=
\lim_{k\to \infty}\| C_{\phi}(F_k)\|_{\mathscr{B}_{\omega_{2}}(\mathbb{B}^n)}=0.
\]
Consequently,  (\ref{composition-Bloch-compact}) follows from $$\lim_{k\to \infty}r_k^{k-1}=\gamma>0.$$

Conversely, assume that (\ref{composition-Bloch-compact}) holds.
Suppose that $\{f_{j}=h_{j}+\overline{g_{j}}\}$ is a sequence of
$\mathscr{B}_{\omega_1}(\mathbb{D})$ such that
$\|f_{j}\|_{\mathscr{B}_{\omega_1}(\mathbb{D})}\leq1$, where $h_{j}$ and
$g_{j}$ are analytic in $\mathbb{D}$ with $g_{j}(0)=0$.
Then there
are subsequences of $\{ h_{j}\}$ and $\{ g_{j}\}$ that converge uniformly on compact
subsets of $\mathbb{D}$ to holomorphic functions $h$ ang $g$, respectively.
Without loss of generality, we assume that the
sequences $\{ h_{j}\}$ and $\{ g_{j}\}$ themselves converge to $h$ and $g$, respectively.
Then the
sequence $f_{j}=h_{j}+\overline{g_{j}}$ itself converges uniformly on compact
subsets of $\mathbb{D}$
to the harmonic function
$f=h+\overline{g}$ with $g(0)=0$.
Consequently,
$$1\geq\lim_{j \rightarrow\infty}\left\{|f_{j}(0)|+\mathscr{B}_{\omega_1}^{f_{j}}(w)\right\}=|f(0)|+\mathscr{B}_{\omega_1}^{f}(w),$$
which gives that $f\in\mathscr{B}_{\omega_1}(\mathbb{D})$ with
$\|f\|_{\mathscr{B}_{\omega_1}(\mathbb{D})}\leq 1$.
Let $\varepsilon>0$ be arbitrarily fixed.
Since  (\ref{composition-Bloch-compact}) holds,
there exists a positive integer $N>m_0$ such that, for
$k\in\{N, N+1,\ldots\}$,
\[
\frac{4}{\gamma}\frac{\omega_1(r_k)}{k}\|\phi^{k}\|_{\mathscr{B}_{\omega_{2}}(\mathbb{B}^n)}<\varepsilon,
\]
which implies that, for $j\geq 1$,
\be\label{ef-1-2}
J_1(f_{j}-f, N)\leq \varepsilon.
\ee
Since the sequence $\{ \Lambda_{f_{j}-f}\}$ converges to $0$ locally uniformly in $\mathbb{D}$,
there exists an integer $j_0$ such that, for $j\geq j_0$,
\be\label{ef-1-3}
J_2(f_{j}-f, N)\leq \varepsilon.
\ee
Combining (\ref{ef-1-1}), (\ref{ef-1-2}) and (\ref{ef-1-3}) gives that, for $j\geq j_0$,
$$
\|C_{\phi}(f_{j}-f)\|_{\mathscr{B}_{\omega_2}(\mathbb{B}^n),s}\leq \varepsilon,
$$
which implies that $C_{\phi}$ is compact.
The proof of this theorem is finished.
\qed

\subsection{The proof of Theorem \ref{Thm-4}}
For  $\alpha,\beta\in(0,1)$, let $f\in\mathscr{PH}(\mathbb{B}^n)\cap\mathscr{L}_{\varpi}(\mathbb{B}^n)$, where $\varpi=\alpha~\mbox{or}~\beta$.  Since $\mathbb{B}^n$
is a simply connected domain, we see that $f$ admits the canonical
decomposition $f=f_{1}+\overline{f_{2}}$, where $f_{1}$ and $f_{2}$
are analytic in $\mathbb{B}^n$ with $f_{2}(0)=0$. Next, we prove $f_{1},~f_{2}\in\mathscr{L}_{\varpi}(\mathbb{B}^n)$.
Let $f=f_{1}+\overline{f_{2}}=u+iv$, where $f_{1}=u_{1}+iv_{1}$ and $f_{2}=u_{2}+iv_{2}$. Then $f\in\mathscr{L}_{\varpi}(\mathbb{B}^n)$ if and only if $u,~v\in\mathscr{L}_{\varpi}(\mathbb{B}^n)$. Let $F=f_{1}+f_{2}$ and $\widetilde{v}=\mbox{Im}(F)$, where ``$\mbox{Im}$" is the imaginary part of a complex number.
Then $iF$ is holomorphic in $\mathbb{B}^n$.
Since $\mbox{Re}(iF)=-\widetilde{v}$,
where ``$\mbox{Re}$" is the real part of a complex number,
and $\mbox{Im}(iF)=\mbox{Im}(if)=u$,
$u$ and $\widetilde{v}$ satisfy
$\mbox{div}(\nabla u)=\mbox{div}(\nabla \widetilde{v})=0$
and
$|\nabla u|=|\nabla \widetilde{v}|$.
Then
by \cite[Corollary 3.11]{N-1991}, we see that
\be\label{eq-gj-0.1}\|\mbox{Re}(iF)\|_{\mathscr{L}_{\varpi}(\mathbb{B}^n),s}=\|\widetilde{v}\|_{\mathscr{L}_{\varpi}(\mathbb{B}^n),s}\leq\,
C\|\mbox{Im}(iF)\|_{\mathscr{L}_{\varpi}(\mathbb{B}^n),s}=C\|u\|_{\mathscr{L}_{\varpi}(\mathbb{B}^n),s}.\ee
By (\ref{eq-gj-0.1}), we have
\be\label{eq-gj-0.2}v_{1}=\frac{v+\widetilde{v}}{2}\in\mathscr{L}_{\varpi}(\mathbb{B}^n)\ee
and
\be\label{eq-gj-0.3}v_{2}=\frac{\widetilde{v}-v}{2}\in\mathscr{L}_{\varpi}(\mathbb{B}^n).\ee
Applying \cite[Corollary 3.11]{N-1991} to $f_{1}$ and $f_{2}$ again, we see that
\be\label{eq-gj-0.4}u_{1}\in\mathscr{L}_{\varpi}(\mathbb{B}^n)\ee
and
\be\label{eq-gj-0.5}u_{2}\in\mathscr{L}_{\varpi}(\mathbb{B}^n).\ee
Consequently, $f_{1}\in\mathscr{L}_{\varpi}(\mathbb{B}^n)$ follows from $(\ref{eq-gj-0.2})$ and $(\ref{eq-gj-0.4})$, and
$f_{2}\in\mathscr{L}_{\varpi}(\mathbb{B}^n)$ follows from $(\ref{eq-gj-0.3})$ and $(\ref{eq-gj-0.5})$.
It follows from \cite[Section 6.4]{Ru-1} that $f_{1},~f_{2}\in\mathscr{B}_{\omega}(\mathbb{B}^n)$, where
$\omega(t)=1/(1-t)^{1-\varpi}$ for $t\in[0,1)$. Hence $f\in\mathscr{PH}(\mathbb{B}^n)\cap\mathscr{L}_{\varpi}(\mathbb{B}^n)$ is equivalent
to $f\in\mathscr{B}_{\omega}(\mathbb{B}^n)$. Hence we only need to show that $C_{\phi}:~\mathscr{B}_{\omega_{1}}(\mathbb{D})\rightarrow\mathscr{B}_{\omega_{2}}(\mathbb{B}^n)$
 is bounded if and only if $$\sup_{k\geq1}\left\{k^{-\alpha}\|\phi^{k}\|_{\mathscr{B}_{\omega_{2}}(\mathbb{B}^n)}\right\}<\infty,$$
and $C_{\phi}$ is compact if and only if
$$\lim_{k\to \infty}\left\{k^{-\alpha}\|\phi^{k}\|_{\mathscr{B}_{\omega_{2}}(\mathbb{B}^n)}\right\}=0,$$
 where $\omega_{1}(t)=1/(1-t)^{1-\alpha}$ and $\omega_{2}(t)=1/(1-t)^{1-\beta}$ for $t\in[0,1)$. By taking $r_{k}=(k-1)/(k-\alpha)$ in Theorem \ref{Thm-Bloch},
  we can obtain the desired result. The proof of this theorem is finished.
\qed
\\

{\bf Data Availability} Our manuscript has no associated data.\\

{\bf Conflict of interest} The authors declare that they have no conflict of interest.

\section{Acknowledgments}
The authors would like to thank the referee for many valuable comments.
The research of the first author was partly supported by the Hunan Provincial Natural Science Foundation of China (grant no. 2022JJ10001), the National Science
Foundation of China (grant no. 12071116), 
 the Key Projects of Hunan Provincial Department of Education (grant no. 21A0429);
 the Double First-Class University Project of Hunan Province
(Xiangjiaotong [2018]469),  the Science and Technology Plan Project of Hunan
Province (2016TP1020), and the Discipline Special Research Projects of Hengyang Normal University (XKZX21002);
 The research of the second author was partly supported by
JSPS KAKENHI Grant
Number JP22K03363.


\begin{thebibliography}{99}
\bibitem{AD}  E. Abakumov and  E. Doubtsov,
Reverse estimates in growth spaces,
\textit{Math. Z.} {\bf 271} (2012), 399--413.




\bibitem{AB}  P. Ahern and J. Bruna,
Maximal and area integral characterization of Hardy-Sobolev spaces in the unit ball of $\mathbb{C}^{n}$,
\textit{Rev. Mat. Iberoam.} {\bf 4} (1988), 123--153.



\bibitem{ABR-2001}  S. Axler, P. Bourdon and W. Ramey,   \textit{ Harmonic function theorem}, Springer-Verlag, New York, Inc, 2001.










\bibitem{CH-2022}
S. L. Chen and H. Hamada,
Some sharp Schwarz-Pick type estimates and their applications of harmonic and pluriharmonic functions,
 \textit{J. Funct. Anal.} {\bf 282} (2022), 109254.

\bibitem{CH-2023}
S. L. Chen and H. Hamada,
On (Fej\'er-)Riesz type inequalities, Hardy-Littlewood type theorems and smooth moduli, \textit{Math. Z.} (2023),
https://doi.org/10.1007/s00209-023-03392-6.


 \bibitem{CHPV}S. L. Chen, H. Hamada, S. Ponnusamy, R. Vijayakumar, Schwarz type lemmas and their applications in Banach spaces,
 \textit{J.  Anal. Math.} (2023), DOI 10.1007/s11854-023-0293-0.


 \bibitem{CHZ2022MZ}
S. L. Chen, H. Hamada and J. -F. Zhu, Composition operators on Bloch
and Hardy type spaces, \textit{Math. Z.} {\bf 301} (2022), 3939--3957.

\bibitem{CPR}  S. L. Chen,  S. Ponnusamy, and A. Rasila,
On characterizations of Bloch-type, Hardy-type, and Lipschitz-type spaces,
\textit{Math. Z.} {\bf 279} (2015), 163--183.






 \bibitem{Du}  P. Duren,
{\it Harmonic mappings in the plane,} Cambridge Univ. Press, 2004.

 \bibitem{DHK-2011}  P. Duren, H. Hamada and G. Kohr, Two-point distortion theorems for harmonic and pluriharmonic mappings,
 \textit{Trans. Amer. Math. Soc.}
{\bf 363} (2011), 6197--6218.

\bibitem{Dy1} K. M. Dyakonov,
Equivalent norms on Lipschitz-type spaces of holomorphic functions,
\textit{Acta Math.} {\bf 178} (1997), 143--167.







\bibitem{Fr} H. Frazer, On the moduli of regular functions, \textit{Proc. London Math. Soc.} {\bf 36} (1934), 532--546.



\bibitem{GZ}  Y. T. Guo and X. J. Zhang,
Composition operators from normal weight general function spaces to Bloch type spaces, submitted.









\bibitem{HL32}
G. H. Hardy and J. E. Littlewood,
Some properties of fractional integrals II,
\textit{Math. Z.}
{\bf 34} (1932), 403--439.





\bibitem{HO}  T. Hosokawa and S. Ohno,
Differences of weighted composition operators acting from Bloch space to $H^{\infty}$, \textit{Trans. Amer. Math. Soc.} {\bf 363} (2011), 5321--5340.



\bibitem{Iz}  A. J. Izzo,
Uniform algebras generated by holomorphic and pluriharmonic
functions, \textit{Trans. Amer. Math. Soc.} {\bf 339} (1993),
835--847.







\bibitem{Ko}  A. Kor\'anyi and S. Vagi, Singular integrals in homogeneous spaces and some problems of classical analysis,
\textit{Ann. Scuola Normale Superiore Pisa} {\bf 25} (1971), 575--648.

\bibitem{KL}  S. G. Krantz and S. Y. Li, Area integral characterizations of functions in Hardy spaces on domains in $\mathbb{C}^{n}$,
 \textit{Complex Variables} {\bf 32} (1997), 373--399.

\bibitem{Kw}  E. G. Kwon, Hyperbolic mean growth of bounded holomorphic functions in the ball,
 \textit{Trans. Amer. Math. Soc.} {\bf 355} (2003), 1269--1294.








\bibitem{MM} K. Madigan and A. Matheson,
Compact composition operators on the Bloch space, \textit{Trans.
Amer. Math. Soc.} {\bf 347} (1995), 2679--2687.






\bibitem{Mo-99} A. Montes-Rodr\'iguez,
The essential norm of a composition operator on Bloch spaces, \textit{Pacific J.
Math.} {\bf 188} (1999), 339--351.

\bibitem{Mo-20} A. Montes-Rodr\'iguez,
Weighted composition operators on weighted Banach spaces of analytic functions, \textit{J. London. Math. Soc.} {\bf 61} (2000), 872--884.

\bibitem{N-1991} C. A. Nolder,
Hardy-Littlewood theorems for solutions of elliptic equations in divergence form,
\textit{Indiana Univ. Math. J.} {\bf 40}(1991), 149--160.


\bibitem{P}  M. Pavlovi\'c,
On Dyakonov's paper Equivalent norms on Lipschitz-type spaces of
holomorphic functions, \textit{Acta Math.} {\bf 183} (1999),
141--143.

\bibitem{Pav-2008} M. Pavlovi\'c,
Derivative-free characterizations of bounded composition operators between Lipschitz spaces,
\textit{Math. Z.} {\bf 258} (2008), 81--86.


\bibitem{Pav} M. Pavlovi\'c,
On the Littlewood-Paley $g$-function and Calder\'on's area theorem,
\textit{Expo. Math.} {\bf 31} (2013), 169--195.




\bibitem{PR-2013} J. A. Pel\'aez and J. R\"atty\"a,
Generalized Hilbert operators on weighted Bergman spaces,
\textit{Adv. Math.} {\bf 240} (2013), 227--267.



\bibitem{PX}  F. P\'erez-Gonz\'alez and J. Xiao, Bloch-Hardy pullbacks,
 \textit{Acta Sci. Math. $($Szeged$)$} {\bf 67} (2001), 709--718.



\bibitem{Ra}  W. Ramey,
Local boundary behavior of pluriharmonic functions along curves,
\textit{Am.  J. Math.} {\bf
108} (1986), 175--191.

\bibitem{RU} W. Ramey and D. Ullrich, The pointwise Fatou theorem
and its converse for positive pluriharmonic functions, \textit{Duke Math. J.} {\bf 49} (1982), 655--675.






\bibitem{Ru-1}  W. Rudin,
\textit{Function theory in $\mathbb{C}^{n}$},
New York: Springer-Verlag, 1980.

\bibitem{Sha} J. H. Shapiro, The essential norm of a composition operator, \textit{Ann. Math.} {\bf 125} (1987), 375--404.


\bibitem{Stei} E. Stein,
Some problems in harmonic analysis, \textit{Proceedings of symposium in pure mathematics,}
{\bf 35} (1979), 3--19.




\bibitem{Vl}  V. S.  Vladimirov,
\textit{Methods of the Theory of Functions of Several Complex
Variables,} (in Russian), M. I. T. Press, Cambridge, Mass., 1966.

\bibitem{WZZ} H. Wulan, D. Zheng and K. Zhu, Compact composition operators on BMOA and the Bloch space, \textit{Proc. Am. Math. Soc.} {\bf 137} (2009), 3861--3868.


\bibitem{Zhao} R. H. Zhao, Essential norms of composition operators between Bloch type spaces, \textit{Proc. Am. Math. Soc.} {\bf 138} (2010), 2537--2546.

\bibitem{Z2}  K. Zhu,
\textit{Operator theory in function spaces},
Monographs and Textbooks in Pure and Applied Mathematics, 139,
Marcel Dekker, Inc., New York, 1990.
\end{thebibliography}
\end{document}